\documentclass[final]{elsarticle}

\makeatletter
\@ifclasswith{elsarticle}{arxiv}{
	\def\ps@pprintTitle{%
		\let\@oddhead\@empty
		\let\@evenhead\@empty
		\def\@oddfoot{\centerline{\thepage}}%
		\let\@evenfoot\@oddfoot}
}{}
\makeatother

\usepackage{amsmath,amssymb,amsthm}
\usepackage{color}
\usepackage{algorithm}
\usepackage{algpseudocode}
\usepackage{mathdots}
\usepackage[utf8x]{inputenc}
\usepackage{pgfplots}
\usepackage{filecontents}
\usepackage{url}

\newcommand{\qt}{\ensuremath{\mathcal{QT}}}

\newtheorem{theorem}{Theorem}[section]
\newtheorem{lemma}[theorem]{Lemma}	

\newtheorem{corollary}[theorem]{Corollary}

\theoremstyle{definition}
\newtheorem{definition}[theorem]{Definition}

\theoremstyle{remark}
\newtheorem{remark}[theorem]{Remark}

\newcommand{\norm}[1]{\lVert #1 \rVert}
\renewcommand{\leq}{\leqslant}
\renewcommand{\geq}{\geqslant}

\begin{document}
	\begin{frontmatter}
		\author[DM,ISTI]{Leonardo Robol\corref{cor1}\fnref{gncs}}
		\ead{leonardo.robol@unipi.it}
		
		\address[DM]{Dipartimento di Matematica, Università di Pisa, Italy.}
		
		\address[ISTI]{Institute of Information Science and
			Technologies ``A. Faedo'', ISTI-CNR, Pisa, Italy.}
		\cortext[cor1]{Corresponding author}
		
		\fntext[gncs]{The research of the author was partially supported by the INdAM/GNCS project ``Metodi di proiezione per equazioni di matrici e sistemi lineari con operatori definiti tramite somme di prodotti di Kronecker, e soluzioni con struttura di rango''. The author is a member of the INdAM Research group GNCS.}
		
		\title{Rational Krylov and ADI iteration
			for infinite size quasi-Toeplitz matrix equations}		
		
		\begin{abstract}
We consider a class of linear matrix equations involving semi-infinite
matrices which have a quasi-Toeplitz structure. These equations arise
in different settings, mostly connected with PDEs or the study
of Markov chains such as random walks on bidimensional lattices. 

We present the theory justifying the existence of the solution
in an appropriate
Banach algebra which is computationally treatable, and we propose 
several methods for computing them. We show how to adapt the ADI
iteration to this particular infinite dimensional setting, and 
how to construct rational Krylov methods. Convergence theory
is discussed, and numerical experiments validate the proposed 
approaches. 
		\end{abstract}
	
	\end{frontmatter}
	
	\section{Introduction}
	
	In this work we are concerned with the numerical solution of linear matrix
	equations of Stein, Lyapunov and Sylvester type, i.e., 
	\begin{equation} \label{eq:lyapsylv}
	  AXB + X + C = 0, \qquad 
	  AX + XA^T + C = 0, \qquad 
	  AX + XB + C = 0, 
	\end{equation}
	where $A,B,C,X$ are infinite dimensional. We restrict our attention to the
	case 
	where $A,B,C$ are semi-infinite quasi-Toeplitz matrices of the form
	\[
	  M = T(m(z))+ E_m, \qquad 
	  T(m(z)) := \begin{bmatrix}
	    m_0 & m_1 & m_2 &  \cdots  \\
	    m_{-1} & m_0 &  \ddots & \ddots \\
	    m_{-2} &  \ddots & \ddots\\
	    \vdots & \ddots 
	  \end{bmatrix},
	\]
	where $E_m$ is a compact operator over $\ell^p$ and $T(m(z))$ is the
	Toeplitz part associated with a symbol $m(z) = \sum_{j \in \mathbb Z} m_j z^j$.
	One of the main motivations
	for this investigation is the solution of quadratic 
	matrix equations, which is often required in the analysis
	of quasi-Birth-Death stochastic processes \cite{bini2018semi,bini2005numerical,bini2018quadratic,bini2017efficient,bini2019solving}. 
	In Section~\ref{sec:numerical} we show that problems of this form
	arise also by discretizing 2D PDEs on unbounded domains. 
	
	In the context of quasi-Birth-Death processes, computing
	the steady state probability can be recast into solving a matrix
	equation $A_{1} G^2 + A_0 G + A_{-1} - G = 0$ \cite{bini2005numerical}. 
	Newton's iteration yields, at every step,
	a Stein equation which can be treated with the
	methods presented in this work. 
	
	Under suitable assumptions, that are further discussed in 
	Section~\ref{sec:preliminaries}, one can approximate matrices in this class 
	at arbitrary accuracy with a finite number of parameters; this
	is achieved by storing the non-negligible coefficients of the symbol
	and a compressed low-rank representation of $E_m$. In addition, it can be proved
	that these matrices form a Banach matrix algebra, a natural
	setting for computing
	matrix functions and solving matrix equations.

	We mention that Lyapunov and Sylvester
	equations are encountered in the solution of 
	other kinds of quadratic equations using 
	Newton's method, such as Riccati equations (either symmetric, yielding a Lyapunov equation, or non-symmetric, yielding a Sylvester one). The results that we
	present can be used as a basis to solve such equations involving
	quasi-Toeplitz matrices.

	\subsection{Related work} 
	
	In the last decades, there has been an increasing interest in the solution
	of linear (and nonlinear) matrix equations of large scale, with matrices
	of increasing dimension. Lyapunov and Riccati equations often arise from the
	study of stability properties of continuous 
	linear time invariant dynamical systems (LTI), 
	and their discrete counterparts give rise to Stein and discrete Riccati
	equations. 
	
	After the introduction of efficient methods for the solution of dense 
	problems of small and medium size in the 70s \cite{bartels1972solution,golub1979hessenberg}, with the development
	of high quality linear algebra libraries such as SLICOT \cite{benner1999slicot}, 
	the research focus has moved to the large scale setting, where $\mathcal O(n^3)$
	methods are not practical. 
	
	In this area, the development of efficient numerical schemes is closely 
	related with rational approximation problems, as can be recognized by
	analyzing ADI or rational Galerkin Krylov methods \cite{beckermann2011error}. These methods typically allow to compute
	solutions of linear matrix equations in compressed form in $\mathcal O(n)$ complexity
	both in storage and in floating point operations count
	\cite{simoncini2016computational}. 
	
	More recently, several authors considered the case of
	infinite matrices, with different purposes. Krylov methods
	have been investigated for solving differential equations acting
	on infinite-dimensional spaces \cite{gilles2018continuous,moret2019krylov}
	and in the context of regularization \cite{novati2017some}. Another interesting
	area where these methods are considered is the computation of 
	eigenvalues in nonlinear problems. Among the many techniques 
	for solving these problems, one is to linearize their Taylor series, as done in the infinite Arnoldi method and similar approaches \cite{gaaf2017infinite,van2016rank}; even if the original 
	problems are defined on a finite dimensional space, this linearization yields an eigenvalue problem with matrices of infinite size. 
	In these cases, to
	be able to construct the space using finite memory and complexity, the authors
	restrict their attention to polynomial Krylov methods (and not rational ones). 
	
	In the present work we show that, when considering quasi-Toeplitz matrices, it is possible (under suitable assumptions)
	to make use of Krylov methods of rational type, which are often superior
	to classical Krylov methods for problems that are not well-conditioned
	\cite{simoncini2016computational}, in particular for the
	task of solving matrix equations.

	\section{Linear matrix equations with infinite quasi-Toeplitz matrices}
	
	\subsection{Preliminaries and notation} \label{sec:preliminaries}
	We denote by $\ell^p$ the space of sequences in $\mathbb C$ whose 
	$p$-th power are summable, that is $(x_j)_{j \in \mathbb N} \in \ell^p$
	if and only if $\sum_{j \in \mathbb N} |x_j|^p < \infty$. As usual, 
	this is extended to the case $p = \infty$ by imposing that 
	$(x_j)_{j \in \mathbb N} \in \ell^\infty$ if and only if it is bounded. 
	On these sequences we consider the usual norm defined by 
	$\norm{x}_p := (\sum_{j \in \mathbb N} |x_j|^p)^{\frac 1p}$ when
	$1 \leq p < \infty$, and $\norm{x}_\infty = \max_{j \geq 0}{|x_j|}$. The
	$\ell^p$ spaces are complete Banach spaces for every
	$p \geq 1$, including $p = \infty$. 
	
	Toeplitz matrices are matrices with constant diagonals: 
	\[
	  T := [a_{j-i}]_{i,j\geq 1} = \begin{bmatrix}
	    a_0 & a_1 & a_2 & \cdots \\
	    a_{-1} & a_0 & a_1 & \ddots \\
	    a_{-2} & a_{-1} & \ddots & \ddots \\
	    \vdots & \ddots & \ddots & \ddots 
	  \end{bmatrix}. 
	\]
	When we do not specify their size, we consider them to be
	semi-infinite, i.e., with rows and columns indexed over the
	positive integers. 
	
	The entries of a semi-infinite
	Toeplitz matrix are completely determined from the ones on the first column
	and the first row, so we can associate them with a Laurent series 
	$a(z) = \sum_{j \in \mathbb Z} a_j z^j$, which is called the
	\emph{symbol} of $T$. In particular, we denote the
	Toeplitz matrix with symbol $a(z)$ with 
	$T(a(z))$, or just $T(a)$ when the context is clear, for 
	improved readability. 
	
	Similarly, we define Hankel matrices, which are constant along
	anti-diagonals, and can be completely determined 
	by the coefficient of a Taylor series $f(z) = \sum_{j \geq 0} f_j z^j$: 
	\[
	  H(f(z)) := [f_{i+j-1}]_{i,j \geq 1} = \begin{bmatrix}
		f_1 & f_2 & f_3 & \cdots \\
		f_2 & f_3 & \ddots & \ddots \\
		f_3 & \ddots & \ddots \\
		\vdots & \ddots 
	  \end{bmatrix}
	\]
	Analogously to the Toeplitz case, we often
	use the compact notation $H(f)$, dropping the dependency on $z$.
	
	We are interested in symbols $a(z)$ such that $T(a)$ defines a bounded operator on the spaces $\ell^p$; for this purpose, it is convenient
	to introduce the \emph{Wiener class}. 
	
	\begin{definition}[Wiener algebra]
		A Laurent series $a(z)$ belongs to the \emph{Wiener algebra} $\mathcal W$
		if its coefficients are summable, i.e., if 
		$\sum_{j\in \mathbb Z} |a_j| < \infty$. We denote the sum of the moduli
		of the coefficient as $\norm{a(z)}_{\mathcal W}$, which defines a norm
		that makes the set $\mathcal W$ a Banach algebra. 
	\end{definition}

	We refer the reader to \cite{boeottcher2005spectral} for a detailed
	discussion on the properties of the Wiener algebra. 
	A remarkable connection between Wiener algebra and Toeplitz operators is that, 
	if $a(z) \in \mathcal W$, then $T(a)$ defines a bounded operator on $\ell^p$
	for all $p \in [1, \infty]$. 
	
	\begin{lemma}[\protect{\cite[Proposition~1.1]{boeottcher2005spectral} and 
		\cite[Theorem~1.14]{boeottcher2005spectral}}] \label{lem:bottchernorm}
		Let $T(a)$ be a Toeplitz operator, and $a(z) \in \mathcal W$. Then, for any
		$1 \leq p \leq \infty$, we have that the $\norm{T(a)}_p \leq \norm{a}_{\mathcal W}$. Moreover, if $p \in \{1, \infty \}$, we additionally have $\norm{T(a)}_{\infty} = \norm{T(a)}_{1} = \norm{a}_{\mathcal W}$. 
	\end{lemma}
	
	When operating on Toeplitz matrices, 
	it is often useful to construct Hankel matrices from the Taylor
	series obtained taking either the positive or the negative coefficients
	of a Laurent series $a(z)$; these are denoted by 
	\[
	  a_+(z) := \sum_{j \geq 0} a_jz^j, \qquad 
	  a_-(z) := \sum_{j \geq 0} a_{-j} z^j. 
	\]
	Hankel matrices are a natural object that appears in the study
	of semi-infinite Toeplitz matrices, as shown in the next result. 
	
	\begin{lemma}[\protect{\cite[Proposition~1.3]{boeottcher2005spectral} 
			and \cite[Proposition~1.2]{boeottcher2005spectral}}] \label{lem:toepprod}
		Let $T(a)$ and $T(b)$ two semi-infinite Toeplitz matrices with symbols in the
		Wiener class $\mathcal W$. Then, 
		\[
		  T(a) T(b) = T(c) - H(a_+) H(b_-)
		\]
		Moreover, $H(a_+) H(b_-)$ is a compact operator on $\ell^p$
		for any $p \in [1, \infty]$. 
	\end{lemma}

	A few results concerning the decay of coefficients of Taylor and Laurent
	series will be useful later on. From now on, we denote by 
	$B(z_0, \rho)$ the ball of center $z_0$ and radius $\rho$, and 
	with $\mathbb S^1$ the unit circle, i.e., $\mathbb S^1 = \partial B(0, 1)$. 
	We use the notation $A_\rho$, where
	$\rho > 1$, to denote the annulus
	\[
	  A_\rho := \{ 
	    z \in \mathbb C \ | \ \rho^{-1} \leq |z| \leq \rho 
	  \}, 
	\]
	so that $A_\rho \subseteq A_{\rho'}$ if $\rho \leq \rho'$. An holomorphic
	function defined on $A_\rho$ for any value of $\rho > 1$, admits an 
	expansion as a Laurent series, and the coefficients decay exponentially.
	
	\begin{lemma}[\protect{\cite[Theorem~4.4c]{henrici1974applied}}] \label{lem:henrici}
		Let $a(z) = \sum_{j \in \mathbb Z} a_j z^j$ 
		  be a Laurent series defined on an annulus $A_\rho$, 
		  with $\rho > 1$. Then, for any $r$ that satisfies $\rho^{-1} < r < \rho$, 
		  \[
		    |a_j| \leq 
		      \max_{|z| = r} |a(z)|\cdot r^{-j} 
		  \]
	\end{lemma}

	Lemma~\ref{lem:henrici} guarantees an exponential convergence
	to zero of the coefficients $|a_j|$ as $|j| \to \infty$. This fact
	can be used to estimate how to chop a Laurent series by ensuring a small
	error in the Wiener norm. 
	
	\begin{lemma}
		Let $a(z)$ be a Laurent series in the Wiener class on the closed annulus $\overline{A_\rho}$. Then, the Laurent polynomial $a_k(z) = \sum_{j = -k}^k a_j z^j$ satisfies
		\[
		  \norm{a(z) - a_k(z)}_{\mathcal W} \leq 2 \max_{|z|\in \{ \rho, \rho^{-1} \}} |a(z)| \cdot \frac{\rho^{-k}}{\rho - 1}. 
		\]
	\end{lemma}

\begin{proof}
	Since $\overline{A_\rho}$ is closed, the
	largest annulus where $a(z)$ is holomorphic is $A_{\rho'} \supseteq A_\rho$, 
	for some $\rho' > \rho$. 
	Thanks to Lemma~\ref{lem:henrici} applied choosing $r = \rho$, we can guarantee that 
	\[
	  |a_j| \leq \max_{|z|\in \{ \rho, \rho^{-1} \}} |a(z)| \cdot \rho^{-|j|}. 
	\]
	Then, we can estimate the error encountered when truncating the Laurent series to the
	central $2k+1$ terms by
	\begin{align*}
	  \norm{a(z) - a_k(z)}_\mathcal{W} &= \sum_{|j| \geq k+1} |a_j| 
	  \leq 
	    2 \max_{|z|\in \{ \rho, \rho^{-1} \}} |a(z)| \cdot \sum_{j=k+1}^\infty \rho^{-j} \\
	    &= 
	    2 \max_{|z|\in \{ \rho, \rho^{-1} \}} |a(z)| \frac{\rho^{-k}}{\rho - 1}.
	\end{align*}
\end{proof}

	As mentioned in the introduction, our interest is in matrices which 
	are Toeplitz plus a compact correction. These can be defined as follows:
	
	\begin{definition}
		The class $\qt_p$ of bounded linear operators from $\ell^p$ into itself, 
		for $1 \leq p \leq \infty$, 
		is defined as follows:
		\[
		  \qt_p := \{ 
		      A = T(a) + E_a, \ | \
		      a(z) \in \mathcal W, \
		      E_a \in K(\ell^p)
		  \}, 
		\]
		where $K(\ell^p)$ is the set of compact operators on $\ell^p$. 
	\end{definition}

	We note that the operators in $\qt_p$ are indeed bounded, in view
	of Lemma~\ref{lem:bottchernorm}. They form a Banach algebra with the
	induced $\ell^p$ norm \cite{bini2018quasi}. In particular, it is possible
	to define matrix functions (through the holomorphic functional calculus \cite{bini2017functions}), 
	and solve linear and quadratic matrix equations (see for instance \cite{bini2018quadratic}) in the class. 
	
	\begin{remark}
	It is often
	useful to employ a slightly different norm, defined by
	$\norm{A} = \norm{T(a)}_\infty + \gamma \norm{E_a}_p$. The class is
	a Banach algebra with respect to this norm 
	provided $\gamma$ is chosen appropriately, and this choice
	simplifies the truncation in the numerical approximation 
	of $A$ (see \cite{bini2018quasi}). For simplicity, we assume to be working 
	with $\ell^p$ norms, but all the results in this paper are easily generalizable to these more general norms.
\end{remark}
	
%
%

	\subsection{Existence of solutions to Sylvester equations} \label{sec:existence}
	
	It is well-known that, for $A,B$ operators on
	a finite dimensional space, the matrix equation
	$AX + XB + C = 0$ has a unique solution if and only if
	the tensorized operator $B^T \otimes I + I \otimes A$ ---
	where $\otimes$ denotes the Kronecker product --- is invertible, 
	which holds if and only if $A$ and $-B$ 
	have disjoint spectra. 
	
	The same result holds when considering $A$ and $B$ as
	operators on a Banach space. However, in this greater generality, 
	the definitions of spectrum, eigenvalues, and eigenvectors are slightly
	more delicate, so to avoid any ambiguity we briefly recall them here. 
	
	\begin{definition}
		Let $A$ be a bounded
		operator defined on a complex Banach space. Then,
		we indicate the set
		\[
		  \sigma(A) := \{ z \in \mathbb C \ | \ zI - A \text{ is not invertible} \}, 
		\]
		as the \emph{spectrum of $A$}. We say that 
		$\lambda$ is an eigenvalue of $A$ is $Ax = \lambda x$ for some 
		$x\neq 0$, which is called an \emph{eigenvector relative to $\lambda$}. We denote by $\Lambda(A)$ the set of eigenvalues 
		of an operator $A$. 
	\end{definition}

	Clearly, $\Lambda(A) \subseteq \sigma(A)$, 
	since if $\lambda$ is an eigenvalue then $\lambda I - A$ is not
	injective. On the other hand, in contrast with the finite dimensional
	case, $\sigma(A)$ is generally a larger set. For instance, consider
	the shifting operator $Z$ on $\ell^p$ defined as follows
	\[
	  Z = \begin{bmatrix}
	    0 & 0 & \cdots \\
	    1 & \\
	    & 1 \\
	    & & \ddots 
	  \end{bmatrix}
	\]
	The operator $Z$ is not surjective, since its image is 
	the subset of sequences $\ell^p$ with the first element
	equal to $0$. On the other hand, $0$ is not an eigenvalue, 
	since $Z$ is injective. So $0$ belongs to the spectrum, 
	but not to $\Lambda(Z)$. Let us recall a few known results
	on the existence of solutions in infinite-dimensional Banach spaces.
	
	
	\begin{theorem}[\protect{\cite[Sylvester--Rosenblum]{bhatia1997and}}]
		\label{thm:existence}
		The Sylvester equation $AX + XB + C = 0$, where
		$A,B,C$ are operators on a complex Banach space,
		has a unique solution $X$ 
		if and only if the spectra of $A$ and $B$ satisfy
		$\sigma(A) \cap \sigma(-B) = \emptyset$.
	\end{theorem}
	
	Under a separation assumption on the spectra of $A$ and $B$, one can explicitly express the solution $X$ in integral form. The following result is 
	often found in the context of finite dimensional linear matrix equations \cite{lancaster1970explicit}, but holds 
	for operators as well \cite[Theorem~2.1]{rosenblum1956operator}. In the context of operators, we need to make some additional assumptions
	on the separation of the spectra. The result holds in any Banach
	algebra, but here is stated for matrices in $\qt_p$. 
	
	\begin{theorem}[\protect{\cite[Theorem~2.1]{rosenblum1956operator}}]
		\label{thm:integralform}
		Let $A,B,C \in \qt_p$, and $\Gamma$ be a union of a finite
		number of Jordan curves in the complex plane such that 
		$(zI - A)^{-1}$ is holomorphic for $z$ inside the components of $\Gamma$, and $(zI + B)^{-1}$ 
		is holomorphic outside. Then, the matrix $X$ defined as
		\[
		  X := -\frac{1}{2\pi i} \int_{\Gamma} (zI - A)^{-1} C (zI + B)^{-1}\ dz
		\]
		belongs to the $\qt_p$ algebra, and solves $AX + XB + C = 0$. 
	\end{theorem}

    When the constant term in the equation is compact, 
    this property is reflected in the solution $X$ as well. This will 
    be an important step for providing a constructive 
    formulation for the solution of such equations in 
    the $\qt_p$ algebra. 
    
    \begin{corollary}
    	\label{cor:compact}
    	Let $AX + XB + C = 0$ a Sylvester equation as in Theorem~\ref{thm:integralform}, and assume that $C$ is 
    	compact. Then, the unique solution $X$ is compact as well. 
    \end{corollary}

\begin{proof}
	In view of Theorem~\ref{thm:integralform}, we can write the
	solution $X$ as 
	\[
	  X = \int_{\Gamma} S(z)\ dz, \qquad 
	  S(z) := (zI - A)^{-1} C (zI + B)^{-1}. 
	\]
	Clearly, $S(z)$ is compact for every $z$. The integral
	is defined through a limit operation of a sequence of 
	compact operators (obtained as combinations of evaluations of $S(z)$), so the closedness of the set compact operators implies
	that $X$ is compact. 
\end{proof}

	As discussed
	later in Section~\ref{sec:toeplitzsplitting}, this result provides the building block
	to develop an algorithm for the numerical solution of the Sylvester and 
	Lyapunov equation. In particular, since the
	solution belongs to $\qt_p$, it can be decomposed as 
	$X = T(x(z)) + E_x$ as Toeplitz-plus-correction, and the two 
	parts can be computed separately. 
	
	\begin{lemma} \label{lem:splitsolutionqt}
		Let $A,B,C \in \qt_p$, and
		 $X$ be the solution to $AX + XB + C = 0$, in the hypotheses of 
		Theorem~\ref{thm:existence}. Then, if 
		\[
		  A = T(a(z)) + E_a, \qquad 
		  B = T(b(z)) + E_b, \qquad 
		  C = T(c(z)) + E_c, 
		\]
		we can decompose $X = T(x(z)) + E_x$ where 
		$x(z) = - c(z) (a(z) + b(z))^{-1}$, and 
		$E_x$ is compact and solves the correction equation
		$A E_x + E_x B + \hat C = 0$, with
		$\hat C := C + A T(x(z)) + T(x(z)) B$. 
	\end{lemma}

	\begin{proof}
		In view of Theorem~\ref{thm:existence}, the solution $X$ is in the space
		$\qt_p$, so it admits a unique decomposition as sum of a Toeplitz
		and a compact correction part $X = T(x(z)) + E_x$. 
		
		Since 
		$AX + XB + C = 0$, and the symbol of sums and multiplication of 
		$\qt$ matrices is the sum and multiplication of their symbols, 
		we need to have $a(z) x(z) + x(z) b(z) + c(z) = 0$. This implies 
		$x(z) = - c(z) (a(z) + b(z))^{-1}$, proving the first claim. 
		
		Then, the fact that $E_x$ solves the correction equation can be 
		verified by 
		using the decomposition
		$X = T(x) + E_x$ and 
		writing
		\[
		  A E_x + E_x B + \underbrace{C+ AT(x) + T(x) B}_{= \hat C} = 0.
		\]
	\end{proof}
	
	\section{Numerical solution of linear matrix equations in $\qt_p$}
	
	We discuss two related numerical algorithms for
	the solution of Sylvester and Stein
	 equations in the space $\qt_p$, one based on Krylov subspaces, 
	 and the ADI iteration; these
	methods are generally well suited to solve linear matrix 
	equations when the solution has good low-rank approximability properties \cite{simoncini2016computational}. In
	our context, this is far from being true: the solution $X$ will generally
	have a non-zero Toeplitz part, which is non-compact. 
	
	The key observation to make
	both approaches feasible is to separate the problem of approximating the
	Toeplitz part and the compact 
	correction, and only use the Krylov (or ADI) method
	for the latter. This section is structured as follows:
	
	\begin{itemize}
		\item 	First, in Section~\ref{sec:toeplitzsplitting} 
		we discuss a procedure
		that computes the Toeplitz part
		$T(x(z))$ of the solution and constructs another Sylvester
		equation with compact right hand side that can be used to recover
		the correction part $E_x$.
		\item  Second, we discuss two methods
		to solve the latter equation: the ADI iteration in 
		Section~\ref{sec:adi}, and a 
		rational Krylov subspace approach in Section~\ref{sec:galerkin}. 
		\item Finally, we show how the method can be adapted to solve
		Stein equations in Section~\ref{sec:stein}. 
	\end{itemize}
	The method based on Krylov subspaces will turn out to be more 
	efficient and robust, but can only be formulated for $p = 2$, 
	since it requires a scalar product and thus to
	be working in an Hilbert space. 
	
	Convergence 
	will be proven for the ADI iteration in Section~\ref{sec:adi} for relevant configurations of the spectrum and in the more general hypothesis
	$1 \leq p \leq \infty$. 
	Then, we show in Section~\ref{sec:galerkin} that 
	the convergence result can be extended when $p = 2$
	to the Galerkin approach as well, 
	and we discuss why the latter has enhanced robustness properties. 
	
	\subsection{Computing the Toeplitz part}
	
	Lemma~\ref{lem:splitsolutionqt} guarantees that the solution
	$X$ can be decomposed as $X = T(x) + E_x$, where $E_x$ is compact and 
	$T(x)$ is Toeplitz. The symbol $x(z)$ is defined as 
	\[
	  x(z) := -c(z) \cdot (a(z) + b(z))^{-1},
	\]
	and we assume that $a(z), b(z)$, and $c(z)$ are given through
	their Laurent series. For $x(z)$ to be well-defined, 
	we need $a(z) + b(z)$ to not vanish on the unit circle $\mathbb S^1 \subseteq \mathbb C$. If this holds, 
	then $a(z) + b(z)$ is invertible in an annulus $A_\rho$, and thus we have the
	existence of the Laurent series for $x(z)$. This assumption
	is not restrictive since, as 
	shown in the next Lemma, is required for the Sylvester equation to be well-posed. 
	
	\begin{lemma} \label{lem:well-posed}
		Let $A,B,C$ be matrices in $\qt_p$, and assume that the Sylvester
		equation $AX + XB + C = 0$ has a unique solution, that is $\sigma(A) \cap \sigma(-B) = \emptyset$. Then, $a(z) + b(z) \neq 0$ for every $z \in \mathbb S^{1}$. 
	\end{lemma}

\begin{proof}
	The condition $A,B \in \qt_p$ implies the existence of decompositions
	$A = T(a) + E_a$ and $B = T(b) + E_b$. We know that 
	$a(\mathbb S^1) = \sigma_{\mathrm{ess}} (T(a))$ \cite{boeottcher2005spectral}. The essential spectrum
	is invariant for compact perturbations, and is always a subset 
	of the spectrum. Hence, we have 
	$a(\mathbb S^1) = \sigma_{\mathrm{ess}}(A) \subseteq \sigma(A)$, and 
	by the  
	the same considerations we obtain $b(\mathbb S^1) \subseteq \sigma(B)$ as well. 
	
	Assume by contradiction that 
	there exists $z \in \mathbb S^1$ such that $a(z) + b(z) = 0$. This implies 
	that $\lambda = a(z) = -b(z)$ is in $\sigma(A) \cap \sigma(-B)$,
	making the Sylvester equation singular in view of Theorem~\ref{thm:existence}. 
\end{proof}
	
	The function $x(z)$ is defined on the annulus 
	$A_\rho = \{ \rho^{-1} < |z| < \rho \}$ for some $\rho > 1$, and 
	can be expanded in a Laurent series. We aim at approximating the
	expansions using a truncated Laurent series as follows:
	\begin{enumerate}
		\item A positive integer $n$ is selected, and we evaluate the function
		  $x(z) z^{n}$ at the $(2n+1)$-th roots of the unity. This operation
		  can be performed efficiently using the FFT. 
		\item Using again the FFT, we interpolate a polynomial
		$p(z)$ that coincides with $x(z) z^n$ at the selected point. 
		\item We set $x^{(n)}(z) := z^{-n} p(z)$ as approximation to $x(z)$. If the
		  approximation is accurate enough, we stop, otherwise we double the
		  value of $n$ and continue. 
	\end{enumerate}

	This procedure is known under the name of 
	\emph{evaluation-interpolation scheme} \cite[Section~8.5]{bini2005numerical}. 
		\label{sec:toeplitzsplitting} 
	It remains to clarify how we 
	decide if a certain approximation is accurate enough. To this aim, we rely
	on the fact that, since the Laurent series is well defined in
	an annulus, the coefficients decay exponentially. So, having 
	prescribed a certain tolerance $\tau$, and given a partial approximation
	\[
	  x^{(n)}(z) = \sum_{j = -n}^n x_j^{(n)} z^j,
	\]
	we stop the iterations if 
	\begin{equation} \label{eq:stop}
	  \sum_{|j|  > \lceil \frac n2 \rceil} |x_j^{(n)}| < \norm{x^{(n)}(z)}_{\mathcal{W}} \cdot 
	  \tau. 
	\end{equation}
	\begin{algorithm}
	\caption{Evaluation-Interpolation scheme for the computation of $x(z)$.}
	\label{alg:evinterp}
	\begin{algorithmic}[1]
		\Procedure{EvInterp}{$a(z), b(z), c(z)$}
		\State $n \gets 4$ \Comment{Set the starting degree for the approximant}
		\State $x \gets 0$
		\While{$\|x(1:\frac n2,\frac{3n}2:2n)\|_1 >= \|x\|_1 \cdot \tau$}
		\State $n \gets 2n$
		\For{$j = 1, \ldots, 2n+1$} 
		\State $v_j \gets -\xi_{2n+1}^{n} c(\xi_{2n+1}^j) / (a(\xi_{2n+1}^j) + b(\xi_{2n+1}^j))$
		\EndFor
		\State $x \gets \Call{IFFT}{v}$.
		\EndWhile
		\EndProcedure
	\end{algorithmic}
\end{algorithm}	

\begin{remark}
	In practice, it is convenient to use \eqref{eq:stop} as stopping criterion, 
	even if it cannot guarantee theoretically that the tails are small, and could 
	be fooled into a premature stopping. We note, however, that if one knows an
	$A_\rho$ where $a(z) + b(z)$ is invertible and $c(z)$ well-defined, 
	then the number of non-negligible coefficients might be bounded using
	Lemma~\ref{lem:henrici}. 
\end{remark}
	
	The resulting procedure is summarized in Algorithm~\ref{alg:evinterp}. 	
	Once $x(z)$ is known, we can make use of the following result to compute
	the correction $E_x$. In the algorithm, the for loop can be replaced
	with an FFT, reducing the cost to $\mathcal O(n\log n)$. 
	
	\begin{lemma} \label{lem:Ex}
		Consider the Sylvester equation $AX + XB + C$ with
		$A = T(a) + E_a$, $B = T(b) + E_b$, $C = T(c) + E_c$. If we set $X = T(x) + E_x$, and 
		$x(z) = -c(z) \cdot (a(z)+b(z))^{-1}$, then $E_x$ solves the 
		Sylvester equation
		\[
		  AE_x + E_x B + \hat C = 0, 
		\]
		where \[
		  \hat C:= 
		  E_c + E_a T(x) + T(x) E_b - H(a_+) H(x_-) - H(x_+) H(b_-). 
		\]
		In particular, $\hat C$ is a compact operator.
	\end{lemma}

\begin{proof}
	By Lemma~\ref{lem:splitsolutionqt} we know that $E_x$ solves the equation
	$AE_x + E_x B + \hat C = 0$, with $\hat C = C + AT(x) + T(x) B$. Then,
	using the relation $T(a) T(b) = T(ab) - H(a_+) H(b_-)$ 
	we obtain
	\[
	  \hat C = T(c + ax + xb) + E_c + E_aT(x) + T(x)E_b - H(a_+)H(x_-) - 
	  H(x_+)H(b_-). 
	\]
	By construction, we have $c(z) + a(z) x(z) + x(z) b(z) \equiv 0$, 
	and hence we have proved the formula for $\hat C$. Note that $\hat C$
	is defined as a linear combination of compact operators, and is therefore
	compact as well. 
\end{proof}

	\subsection{An ADI iteration}
	\label{sec:adi}
	
	The ADI iteration has received considerable attention in the last decades
	as an iterative method for the solution of linear matrix equations, in 
	particular because of the connection between the error representation
	and rational approximation problems. We refer the reader to \cite{simoncini2016computational} and the references therein for a complete
	overview. 
	
	The ADI iteration for solving a Sylvester equation
	\begin{equation} \label{eq:sylv}
		AX + XB + C = 0
	\end{equation}
	 is defined by setting 
	$X^{(0)} := 0$ and the subsequent iterates as:
	\begin{align*}
	X^{(j + \frac 12)} &:= -(A - \beta_{j+1})^{-1} X^{(j)} (B + \beta_{j+1}) - (A - \beta_{j+1})^{-1} C \\
	X^{(j + 1)} &:= - (A - \alpha_{j+1}) X^{(j+\frac 12)} (B + \alpha_{j+1})^{-1} - C (B + \alpha_{j+1})^{-1}, 
	\end{align*}
	where $\alpha_j, \beta_j$ are the ADI parameters, 
	and we define a half-step iterate $X^{(j+\frac 12)}$ for
	notational convenience. One can obtain $X^{(j+1)}$ as a 
	function of $X^{(j)}$ directly by composing the above formulas. 
	
	When $C$ has rank $k$, the
	iteration can be rephrased so that $X^{(j)}$ is expressed 
	in a factored low-rank form. In particular,
	if $X^{(0)} = 0$ then $X^{(j)}$ has rank bounded by $jk$. We
	refer the reader to \cite{benner2014computing} for further details. 
	
	One of the key properties of ADI is its error representation, 
	which is closely related with the choice of the parameters
	$\alpha_j, \beta_j$. Let us define the family of rational functions
	\[
	r_k(z) := \prod_{j = 1}^k \frac{z - \alpha_j}{z - \beta_j}, 
	\qquad 
	k = 1, 2, \ldots 
	\]
	Then, the error at step $k$ of the ADI method can be written as:
	\begin{equation} \label{eq:adierror}
      X - X^{(k)} = r_k(A) X r_k(-B)^{-1}. 
	\end{equation}
	This gives an indication on how to choose the parameters $\alpha_j, \beta_j$, assuming one is able to solve a rational approximation
	problem and find a rational function of degree $k$ that is small
	when evaluated at $A$, and such that its inverse is small if evaluated at $-B$. In case of 
	normal matrices when working with unitarily invariant norms, this 
	can be recast as a problem on the eigenvalues and is known
	as a \emph{Zolotarev problem}. An explicit solution is difficult
	to find in a general context, but has been given by Zolotarev in 
	1877 for the case of two real intervals $[a, b]$ and $[-b, -a]$
	\cite{zolotarev1877application}. 	
	For more general cases, some
	heuristics are described in \cite{sabino2006solution}. When
	considering the $2$-norm and normal
	matrices $A$ and $B$, \eqref{eq:adierror} also gives the
	following error bound on the relative residual
	\begin{equation} \label{eq:adiresidual}
	  \frac{\norm{AX^{(k)} + X^{(k)} B + C}_2}{\norm{X}_2}
	     \leq (\norm{A}_2 + \norm{B}_2) \cdot \frac{
	     	\max_{\lambda \in \sigma(A)} |r_k(\lambda) |
	     }{
	      	\min_{\lambda \in \sigma(B)} |r_k(-\lambda) |
 		}
	\end{equation}
	If considering Hermitian positive definite matrices $A$ and $B$
	the result by Zolotarev yields (see \cite{beckermann2011error} for 
	a more modern reference):
	\begin{equation} \label{eq:adizolotarev}
	  \frac{\norm{AX^{(k)} + X^{(k)} B + C}_2}{\norm{X}_2}
	\leq 4 (\norm{A}_2 + \norm{B}_2) \rho^{k}, \qquad 
	\rho = e^{-\frac{\pi^2}{\log(4\frac{b}{a})}}. 
	\end{equation}
	
	This method is applicable using the arithmetic of 
	the $\qt_p$ class directly, since it involves the solution
	of linear systems and matrix multiplications by (shifted)
	$\qt_p$ matrices. Efficient implementation of these arithmetic
	operations are provided, e.g., in the \texttt{cqt-toolbox}\footnote{Available at \url{https://github.com/numpi/cqt-toolbox}.} for MATLAB \cite{bini2018quasi}. When starting from $X_0 \equiv 0$, and 
	a low rank (resp. compact) right hand side $C$, the solution will be 
	numerically low-rank (resp. compact), and the implementation of
	arithmetic operations in \texttt{cqt-toolbox} automatically exploits this
	property (for further details, we refer to the recompression
	techniques described in \cite{bini2018quasi}). 
	
	To make the method attractive, a good choice
	for the parameters $\alpha_j, \beta_j$ is needed.
	When these are available, we can provide a bound for the accuracy of 
	$X$ by bounding the infinity norm of $f(A)$ 
	where $f(z)$ is a rational function	$r_k(z)$. Recall
	that we use the notation $B(z_0, \rho)$ to denote the open ball
	of radius $\rho$ centered at $z_0$. 
	
	\begin{lemma} \label{lem:bound}
		Let  $f(z)$ be a holomorphic function on $B(z_0, \rho) \subseteq \mathbb C$, 
		$z_0 \in \mathbb C$, 
		$A$ be a matrix such that $\|A - z_0 I\| \leq \rho' < \rho$,
		with $\|\cdot\|$ being any subordinate norm.
		Then, 
		\[
		\| f(A) \| \leq \frac{\rho}{\rho - \rho'} \max_{z \in B(z_0, \rho)}
		|f(z)|. 
		\]
	\end{lemma}
	
	\begin{proof}
		Note that, without loss of generality, we can assume that $z_0 = 0$. Indeed,
		if that's not the case, we may set $\tilde A := A - z_0 I$, and $\tilde f(z) := f(z + z_0)$, and we have that $\|\tilde A\| \leq \rho' \iff \|A - z_0 I\| \leq \rho'$, and 
		\[
		  \norm{f(A)} = \norm{\tilde f(\tilde A)} \leq \frac{\rho}{\rho - \rho'} \max_{z \in B(0, \rho)} |\tilde f(z)| = \frac{\rho}{\rho - \rho'} \max_{z \in B(z_0, \rho)} |f(z)|. 
		\]
		We consider the Cauchy integral representation of $f(A)$,
		given by
		\[
		f(A) = \frac{1}{2\pi i} \int_{\partial B(0, \rho)} f(z) (zI - A)^{-1}\ dz, 
		\]
		Since $\norm{\cdot}$ is an induced norm we have
		 $\rho(A) \leq \norm{A} \leq \rho'$, so $B(0, \rho)$
		contain the spectrum of $A$. Then, we can bound the integral
		taking the maximum of the norm (or absolute values) of the integrand, multiplied by the length of the integration path, 
		which yields
		\[
		\| f(A) \| \leq \rho \cdot \max_{z \in B(0, \rho)} |f(z)| \max_{z \in \partial B(0, \rho)} \| (zI - A)^{-1} \|, 
		\]
		where we have used that the maximum of the modulus of $f(z)$ 
		on the boundary is equal to the maximum inside the set, in view
		of the maximum modulus principle. Then, we have 
		that for $|z| = \rho$, 
		\[
		(zI - A)^{-1} = z^{-1} (I - z^{-1} A)^{-1} = 
		z^{-1} \sum_{j = 0}^{\infty} (z^{-1} A)^j \implies 
		\|(zI - A)^{-1}\| \leq \frac{1}{\rho - \rho'},
		\]
		where the last equality follows by taking the norms
		and using $\|z^{-1} A\| \leq \rho' / \rho < 1$ for $|z| = \rho$. 
		Combining these estimates yields the desired bound. 
	\end{proof}

	The above estimate allows to give a general statement of ADI convergence.
	
	\begin{theorem} \label{thm:adierr}
		Let $X^{(k)}$ the
		approximation to the solution $X$ of $AX + XB + C = 0$
		obtained after $k$ steps of ADI with parameters
		$\{ \alpha_j \}$ and $\{ \beta_j \}$, $z_A, \rho_A, \rho_A'$ such that 
		$\| A - z_A I \| \leq \rho_A' < \rho_A$, and $z_B, \rho_B, \rho_B'$ such that 
		$\| B - z_B I \| \leq \rho_B' < \rho_B$. 
		Then, if we define the
		rational function $r_k(z) := \prod_{j = 1}^k (z - \alpha_j) / (z - \beta_j)$, 
		\[
		  \norm{X - X^{(k)}} \leq \norm{X}\cdot \frac{\rho_A \rho_B}{(\rho_A- \rho_A')(\rho_B- \rho_B')} \frac{\displaystyle \max_{z \in B(z_A, \rho_A)} |r_k(z)|}{\displaystyle \min_{z \in B(z_B, \rho_B)} |r_k(-z)|}, 
		\]
		where $\norm{\cdot}$ is any induced norm. 
	\end{theorem}

	\begin{proof}
		In view of \eqref{eq:adierror}, writing the error at step $k$ of  ADI and taking norms yields
		\[
		  \norm{X - X^{(k)}} \leq \norm{X} \cdot \norm{r_k(A)^{-1}} \cdot 
		   \norm{r_k(-B)}. 
		\]
		We can bound $ \norm{r_k(A)}$ using Lemma~\ref{lem:bound} with
		$f(z) = r_k(z)$, which yields
		\[
		  \norm{r_k(A)} \leq \frac{\rho_A}{\rho_A - \rho_A'} \max_{z \in B(z_A, \rho_A)} |r_k(z)|. 
		\]
		Similarly, we can bound the term $\norm{r_k(-B)^{-1}}$ setting 
		$f(z) = r_k(-z)^{-1}$:
		\[
		\norm{r_k(-B)^{-1}} \leq \frac{\rho_B}{\rho_B - \rho_B'} \max_{z \in B(z_B, \rho_B)} |r_k(-z)^{-1}| = \frac{\rho_B}{\rho_B - \rho_B'} \frac{1}{\displaystyle \min_{z \in B(z_B, \rho_B)} |r_k(-z)|}. 
		\]
		Combining these two inequalities concludes the proof. 
	\end{proof}
	
%
%
	
	\subsection{Rational Krylov Galerkin approximation}
	\label{sec:galerkin}
	
	A limitation of the ADI approach is that, whereas good choices
	of the parameters yield a fast convergence, the method is not
	robust to perturbation in these values \cite{beckermann2011error,simoncini2016computational}. 
	If the parameters $\alpha_j$ and $\beta_j$ are slightly changed, 
	or if they are 
	not chosen in an optimal way for the
	problem at hand, the convergence can be easily degraded. 
	
	We propose to use a strategy for the solution that partially overcomes this
	limitation -- but at the same time requires the stronger hypotheses
	that the correction $E_x$ to be calculated is not a compact
	operator on $\ell^{p}$ for a generic $p$, 
	but on the Hilbert space $\ell^2$. 
	
	With this additional hypothesis, 
	we can employ a Galerkin approach, where only the
	poles $\beta_j$ are chosen, and the method has the same convergence properties
	as ADI with the best possible
	$\alpha_j$ (in a least square sense) \cite{beckermann2011error}. 
	
	
	In order to introduce the method, we first need to recall
	the definition of a (block) Krylov subspace, denoted
	by $\mathcal K_m(A, U)$:
	\[
	\mathcal K_m(A, U) := \mathrm{span}(
	U, AU, \ldots, A^{(m-1)}U 
	).
	\]
	If $U$ has $k$ columns and rank $k$, we generically expect
	the dimension of $\mathcal K_m(A, U)$ to be $km$, 
	even though in particular cases deflation or breakdown
	might occur. Consider a set of poles $\beta_1, \ldots, \beta_m$, 
	define $q(z) = \prod_{j = 1}^m(z - \beta_j)$, and 
	assume the $\beta_j \not\in \sigma(A)$. Then, 
	the rational Krylov subspace associated with these poles can
	be defined as 
	\[
	\mathcal {RK}_m(A, U, \{ \beta_1, \ldots, \beta_m \}) :=
	q(A)^{-1} \mathcal K_m(A, U). 
	\]
	Formally, we may choose some 
	$\beta_j = \infty$, which by convention means that there
	is a degree deficiency in $q(z)$. 
	
	Independently of the choice for the space, 
	we can formulate the Galerkin projection method as follows. Denote by $W_m$ and $Z_m$ the matrices whose columns form an orthogonal basis of the
	selected
	subspaces; we recover an approximate solution to \eqref{eq:sylv} by solving the
	projected equation
	\[
	(W_m^* A W_m) Y_m + Y_m (Z_m^* B Z_m) + W_m^* C Z_m = 0
	\]
	obtained multiplying \eqref{eq:sylv} by $W_m^*$ on the left and $Z_m$ on the right. This 
	equation has finite dimension, and can therefore be solved by a dense solver such as the
	Bartels-Stewart algorithm \cite{bartels1972solution}. 
	Then, we consider $X_m = W_m Y_m Z_m^*$
	as an approximation to $X$. 
	
	The following result shows that, when $C$ has low-rank, 
	there is a close connection between the Galerkin projection
	method with a set of poles, and the ADI method with the same poles. 
	In the $2$-norm, one can show that the residual of the
	solution recovered by Galerkin is only worse up to a constant
	compared to the one obtained by ADI. Since this holds independently
	of the numerator, this means that Galerkin will match the
	accuracy of ADI with the best possible choice of $\alpha_j$ --- and 
	this makes the method much more robust to variations in the
	parameters (we have to choose only $m$ of them, and not $2m$: the 
	others are automatically determined in a quasi-optimal way). 
	\begin{theorem}
		Assume $C = UV^*$, and let 
		$X^{(k)}_{\text{G}}$ be the solution obtained
		using the Galerkin
		projection method for solving \eqref{eq:sylv} using the rational
		Krylov spaces
		\[
		  \mathcal{RK}(A, U, \{ \beta_1, \ldots, \beta_m \}), \qquad 
		  \mathcal{RK}(B^*, V, \{ -\alpha_1, \ldots, -\alpha_m \}),
		\]
		and 
		$X^{(k)}_{\text{ADI}}$ the solution obtained using ADI with parameters
		$\alpha_j, \beta_j$. Then, 
		\[
		\| AX^{(k)}_{\text{G}} + X^{(k)}_{\text{G}}B + C \|_2 \leq (1 + \kappa_2(A) + \kappa_2(B))
		\| AX^{(k)}_{\text{ADI}} + X^{(k)}_{\text{ADI}}B + C \|_2, 
		\]
		where $\kappa_2(A) := \norm{A}_2 \cdot \norm{A^{-1}}_2$ is the
		condition number in the $\ell^2$-norm. 
	\end{theorem}

\begin{proof}
	The result has been first proven in the Frobenius norm in  \cite[Theorem~2.1]{beckermann2011error}. The extension to the Euclidean
	norm can be found in \cite{massei2019kronecker}, and 
	the latter proof holds unchanged in the Hilbert space $\ell^2$. 
\end{proof}
	
	Using Theorem~\ref{thm:adierr}, we can characterize
	the choice of poles
	that give a fast convergence for ADI, and therefore it is 
	reasonable to make the same choice
	for the Galerkin projection method. 
	
	As reported in the 
	numerical experiments, this choice produces good performances, 
	and the quasi-optimality of the Galerkin projection (that automatically
	optimizes the numerator in Theorem~\ref{thm:adierr})
	often accelerates the convergence. 
	
	\subsection{Solution of Stein equations} \label{sec:stein}
	
	The presented algorithms can be adapted 
	for solving Stein equations of the form 
	\begin{equation} \label{eq:stein}
	MXN + X + C = 0. 
	\end{equation}
	The two problems are closely related, and the Stein equation is
	solvable if and only if $\sigma(M) \cap \sigma(-N)^{-1} = \emptyset$, where
	by $\sigma(N)^{-1}$ we denote the inverses of the elements in the
	spectrum of $N$. 
	
	Indeed, \eqref{eq:stein} has a unique solution
	if and only if the operator $X \mapsto MXN + X$ is invertible, and we have
	$\sigma(X \mapsto MXN) \subseteq \sigma(M)\sigma(N)$ (for operators on 
	a Banach space the inclusion can be obtained following the same proof of existence
	and uniqueness for Sylvester equations in \cite{bhatia1997and}). 
	In particular, the 
	spectrum of the Stein operator is enclosed in $\sigma(M)\sigma(N) + 1$, 
	so it does not contain the zero if and only if $\sigma(M) \cap \sigma(-N)^{-1} = \emptyset$. 
	
	Under the stronger assumption that there exists a disc of radius $\rho < 1$ such
	that $\sigma(M), \sigma(N) \subseteq B(0, \rho)$, 
	we can describe the solution more explicitly. 
	
	\begin{lemma} \label{lem:fixed-point-stein}
		Let $M,N$ operators on a Banach space, such that $\sigma(M), \sigma(N) \subseteq B(0, \rho)$ for some $0 < \rho < 1$. Then, the solution $X$
		to \eqref{eq:stein} is unique and given by 
		\[
		  X = \sum_{j = 0}^\infty (-1)^{j+1} M^j C N^j. 
		\]
		In addition, if $C$ is compact then $X$ is compact as well, and 
		if $C \in \qt_p$ then $X \in \qt_p$. 
	\end{lemma}

	\begin{proof}
		The existence and uniqueness follows by the previous considerations, since
		by construction $\sigma(M) \subseteq B(0, \rho)$, $\sigma(-N)^{-1} \subseteq B(0,\rho^{-1})^C$, and $B(0,\rho) \cap B(0,\rho^{-1})^C = \emptyset$. The explicit expression of the solution $X$ can be
		given by recursively expanding the matrix iteration
		\[
		  X_{k+1} = - C - M X_k N. 
		\]
		To prove that the iteration converges to $X$, note that $X = -C - MXN$ 
		and subtract $X$ 
		from the left hand side and $-C-MXN$ from the right hand side; we obtain
		the error equation
		\[
		  X_{k+1} - X = - M (X_k - X) N, \implies 
		  \norm{X_{k+1} - X} \leq \norm{M}\cdot \norm{N} \cdot \norm{X_k - X}. 
		\]
		Hence, the iteration convergences geometrically to $X$ 
		with rate $\rho^2$ because $\norm{M} \norm{N} < \rho^2$. In addition if $C$ is compact then 
		$Y_{k} := \sum_{j = 0}^k (-1)^{j+1} M^j C N^j$ is compact as well, 
		and we conclude noting that $\norm{X- Y_k} \to 0$ for $k \to \infty$. 
		The set $\qt_p$ is closed, so the second part of the statement follows
		using the same argument. 
	\end{proof}

	The previous Lemma provides an algorithm for
	the computation of $X$. An alternative approach can be given by
	recasting the problem into a Sylvester equation. 
	
	\begin{lemma} \label{lem:cayley}
		Let $M, N$ be matrices with eigenvalues contained in the disc $B(0, \rho) = \{ 
		z \in \mathbb C \ | \ |z| \leq \rho
		\}$, with $\rho < 1$. Then, the matrices 
		\[
		A := (M + I) (I - M)^{-1}, \qquad B := (N + I)^{-1} (I - N)
		\]
		satisfy $\sigma(A) \subseteq - C_\rho$ and $
		\sigma(B) \subseteq C_\rho$, where $C_\rho \subseteq \{ z \in \mathbb C \ | \ \Re(z) > 0\}$ is the disc of center $(\rho^2 + 1) / (1 - \rho^2)$ 
		and radius $2\rho /(1 - \rho^2)$. 
	\end{lemma}
	
	\begin{proof}
		Note that the matrices $A$ and $B$ are obtained by $M$ and $N$ through
		the M\"obius transformation
		\[
		  \mathcal C(z) := \frac{z + 1}{1 - z}
		\]
		by setting $A = \mathcal C(M)$ and $B = \mathcal C(-N)$. Since M\"obius transforms
		map circles into circles (or, more generally, projective lines), one can
		easily verify that:
		\[
		  \mathcal C(B(0, \rho)) = B\left( \frac{\rho^2 + 1}{1 - \rho^2}, \frac{2\rho}{1 - \rho^2} \right) \subseteq \{ z \ |\  \Re(z) > 0 \}. 
		\]
		Indeed, since the image of $B(0, \rho)$ needs to be symmetric
		with respect to the real axis, the formula for the
		center of the resulting disc
		is given by $\frac{1}{2}(\mathcal{C(-\rho)} + \mathcal{C(\rho)})$; 
		analogously, 
		the radius can be computed by $\frac{1}{2}(\mathcal{C(\rho)} - \mathcal{C(-\rho)})$. Since the spectra of $M$ and $N$ 
		are mapped through $\mathcal C$ (in view of the spectral
		mapping theorem), and the inclusions are preserved, this concludes the
		proof. 
	\end{proof}

	The previous result enables the following reformulation of the problem. Instead
	of considering \eqref{eq:stein}, we may consider the Sylvester equation
	$
	  AX + XB + \tilde C = 0, 
	$
	where 
	\[
	  A = (M+I)(I-M)^{-1}, \quad 
	  B = (N+I)^{-1}(I-N), \quad 
	  \tilde C =2 (I - M)^{-1} C (I+N)^{-1}.
	\]
	A direct computation shows that the two equations are equivalent. 
	
	The next result shows that the poles
	for the ADI iteration can be chosen to ensure at least the same speed of
	convergence of the fixed point iteration described in Lemma~\ref{lem:fixed-point-stein}.
	
	\begin{lemma} \label{lem:adiconv}
		Let $X^{(k)}$ be the approximation to the solution of $MXN + X + C = 0$ 
		obtained after $k$ steps of ADI applied to the Sylvester equation
		\eqref{eq:sylv} obtained by remapping the coefficients as described 
		in Lemma~\ref{lem:cayley}, and using the parameters
		$\alpha_j = 1, \beta_j = -1$. Then, for every $\rho$ 
		such that $\|M\|, \|N\| < \rho < 1$, we have
		\[
		\norm{X - X^{(k)}}_p \leq 
		\norm{X}_p \cdot \frac{\rho^{2k+2}}{(\rho - \|M\|)(\rho - \|N\|)}  
		\]
	\end{lemma}

	\begin{proof}
		In view of Theorem~\ref{thm:adierr}, we may consider the
		rational functions $r_k(z) = (z + 1)^k / (z - 1)^k$, 
		and we have the error estimate
		\[
		  \norm{X - X^{(k)}}_p \leq 
		  \norm{X}_p \cdot \frac{\rho^{2}}{(\rho - \|M\|)(\rho - \|N\|)} \cdot 
		    \frac{\max_{z \in C_\rho} |r_k(z)|}{\min_{z \in C_\rho} |r_k(-z)|}
		\]
		It is easy to check that $r_k(z) = (\mathcal C^{-1}(z))^k$. Therefore, 
		$\max_{z \in C_\rho} |r_k(z)| = \max_{z \in B(0, \rho)} |z^k| = \rho^k$. 
		Similarly, a direct computation yields $r_k(-z) =  (\mathcal C^{-1}(z))^k$, 
		and therefore $\min_{z \in C_\rho} |r(-z)| = \rho^{-k}$. Combining
		these identities we have the sought bound. 
	\end{proof}

\begin{remark}
	We note that for the result to be the most effective, it is advantageous
	to scale $M$ and $N$ to have the same norm. This can always be done, 
	choosing a positive scalar $\lambda$ imposing
	$
	  \norm{\lambda M}_p = \norm{\lambda^{-1} N}_p = \sqrt{\norm{M}_p \norm{N}_p}, 
	$ which yields $\lambda = \sqrt{\norm{M}_p^{-1} \norm{N}_p}$. 
\end{remark}

	\begin{remark}
		The remapping of a Stein equation into a Sylvester one is a
		standard tool in the solution of linear matrix equations, see 
		for instance \cite{simoncini2016computational}. In our setting, 
		both the iteration of Lemma~\ref{lem:fixed-point-stein} and 
		the ADI or Galerkin solvers can be efficiently used in the $\qt$ 
		arithmetic. The former is often advantageous for its simplicity, while
		the latter can easily handle more general spectral configurations, 
		and benefits from the Galerkin acceleration obtained by the improved
		robustness with respect to the pole choice. 
	\end{remark}
	
	\section{Rational Krylov methods and $\qt$ matrices}
	\label{sec:rkqt}
	
	The ADI and the fixed point iteration of Lemma~\ref{lem:fixed-point-stein}
	 for the Stein case 
	can be directly performed in $\qt$ arithmetic relying on the
	operations implemented in \texttt{cqt-toolbox}. 
	
	Implementing the Galerkin projection, on the other hand, is more challenging  in an infinite-dimensional settings, since an ad-hoc
	procedure for the representation of the Krylov basis and recurrence
	needs to be developed. 
	
	\subsection{Representing infinite dimensional subspaces}
	\label{sec:storingbasis}
	It 
	is our interest to compute an orthogonal basis of the rational
	Krylov subspace $\mathcal{RK}_m(A, U)$, and this task needs to 
	accomplished considering that $U$ is a vector in an infinite
	dimensional space (but with finite support), and $A$ is an infinite
	quasi-Toeplitz matrix. 
	
	Then, the basis is constructed using the rational Arnoldi process
	\cite{ruhe1984rational}, which requires to compute matrix vector products
	$Ax$ and linear system solutions $(A - \gamma I) x = b$. In the
	$\mathcal QT$ arithmetic, one can represent the vectors $x$ and $b$
	as infinite matrices with zero Toeplitz part, where just the
	first column is non zero. Within this setting, we already have
	the necessary truncation and approximation procedures to compute
	$Ax$ and $(A - \gamma I)^{-1} b$ --- and therefore we can rely
	on the algorithms implemented in \texttt{cqt-toolbox}. For instance, 
	the solution of a linear system can be performed with the
	following instruction:
	\begin{verbatim}
	>> cB = cqt([], [], b); % Zero symbol, only the correction 
	>> cX = A \ cB; % Here A is a CQT matrix
	>> x = correction(cX); % We extract the result
	                       % reading the correction
	\end{verbatim}
	At first sight, it might seem inefficient to store a vector as a matrix
	with only one non-zero column. However, considering how the storage is
	implemented in \texttt{cqt-toolbox}, this is very efficient --- the zero
	part of the matrix is just ignored and not stored at all \cite{bini2018quasi}.
	
	Concerning the reorthogonalization, we employ a 
	modified Gram-Schmidt orthogonalization procedure, that
	requires scalar products and normalization, 
	which can always be easily performed relying on the toolbox. 
	
	The Galerkin rational Krylov solver for linear matrix equation has been
	included in the most recent release of \texttt{cqt-toolbox}, and 
	is available under the name \texttt{rk\_sylv}. Our implementation is 
	similar to the one found in \texttt{rktoolbox} \cite{berljafa2015generalized}, modified to properly handle scalar products 
	and truncations needed in the infinite-dimensional setting. 
	
	\section{Numerical experiments} \label{sec:numerical}
	
	In this section we test our computational framework on two representative
	examples. This confirm our theoretical findings, and also demonstrate the
	differences in efficiency between the different approaches. 
	
	In particular, it is verified that the Galerkin approach is often the most
	robust and practical choice. 
	
	\subsection{A PDE on an unbounded domain}

	We consider the following 2D Poisson problem on the positive
	orthant:
	\begin{equation} \label{eq:pde}
	  \begin{cases}
	    \frac{\partial u(x,y,t)}{\partial t} = \Delta u(x,y,t) + f(x,y,t) & 
	      \text{on } \mathbb R_+ \times \mathbb R_+ \\
	    u(x,y,t) \equiv 0 & \text{on } \mathbb R_+ \times \{0\} \cup \{ 0\} \times \mathbb R_+ \\
	  \end{cases}. 
	\end{equation}
	Given a spatial discretization step 
	$\Delta x$, we consider the finite difference discretization
	obtained using central differences for the second derivative, which
	yields the semi-discretization
	\[
	  \begin{cases}
	  u'(t) = \mathcal M u(t) + f(t), & t > 0\\
	  u(t) \equiv 0 & t = 0, 
	  \end{cases}
	\]
	where with a slight abuse of notation we have denoted by
	$u(t)$ and $f(t)$ the discretized-in-space versions of $u(x,y,t)$ and 
	$f(x,y,t)$. Thanks to the structure of the problem, the matrix $\mathcal M$
	has the special structure 
	$
	  \mathcal M = A^T \otimes I + I \otimes A, 
	$
	where $A$ discretizes the 1D second derivative, and is equal (up to
	a scaling)
	to the Toeplitz tridiagonal matrix with symbol $a(z) = z^{-1} - 2 + z$, . Hence, solving a linear
	system with $\mathcal M$ is equivalent to solving a matrix equation 
	with coefficients equal to $A$. 
	
	We now discretize in time the problem
	with an implicit method; for simplicity, we consider the 
	implicit Euler method with time step $\Delta t$. 
	This yields the following time stepping scheme starting from $u(x,y,0) \equiv 0$:
	\[
	  (\Delta t \mathcal M - \Delta x^2 I) u_{t + 1} = -\Delta x^2 u_t - \Delta x^2
	  \Delta t f_{t+1}, 
	\]
	where $u_t$ is the discretized in space solution at time $t$, and 
	$f_t$ the evaluation of $f(x,y,t)$ on the grid at the same instant in time. 
	Exploiting the Kronecker structure of $\mathcal M$, this linear system
	can be rephrased (up to rearranging the entries of $u_t$ and $f_{t+1}$ in the
	matrices $U_t$ and $F_{t+1}$) as the linear matrix equation:
	\[
	  \left( \Delta t A - \frac{\Delta x^2}{2} I \right) U_{t+1} + 
	  U_{t+1} \left( \Delta t A - \frac{\Delta x^2}{2} I \right)
	  = -\Delta x^2 U_t - \Delta x^2 \Delta t \cdot F_{t+1}, 
	\]
	We consider the source function
	\[
	  f(x,y) = \frac{1}{10} e^{-(x - y)^2} + e^{-(x+y)^2},
	\] whose sampling on the grid is the sum
	of a Toeplitz and an Hankel matrix, and is therefore in the $\qt$ class. 
	
	We note that $A$ is negative semi-definite, and shifting it with the 
	identity makes it negative definite. This kind of reformulation of linear
	systems arising from PDEs into linear matrix equation is well-known, 
	see \cite{kressner2019low,massei2018fast,massei2018solving,palitta2016matrix} 
	for some recent works where the problem is studied and examples with
	(finite) Toeplitz matrices are given. 
	
	The matrix $F$ sampling $f(x,y)$ 
	can be constructed as follows: we compute all the coefficients
	of the symbol of $e^{-(x - y)^2}$ (as a Toeplitz matrix) and of 
	$e^{(x+y)^2}$ (as a Hankel matrix). Then, we construct the Toeplitz
	part using the constructors in \texttt{cqt-toolbox}, and we use the
	function \texttt{cqt('hankel', ...)} to construct a \texttt{cqt} representation
	of the Hankel part. In practice, the construction is done as follows:
	\begin{verbatim}
% Right hand side
fm = exp(-(0 : dx : log(1/sqrt(eps))).^2);
F = cqt('hankel', fm) + 0.1 * cqt(fm,fm);
	
% Laplacian with Dirichlet boundary conditions
A = cqt([-2 1], [-2 1]);
	
% Perform some time stepping
X = cell(1, timesteps + 1);
X{1} = cqt([], []);
		
for j = 1 : timesteps
    M = dx^2/2 * cqt(1, 1) - dt * A;
    C = dx^2 * X{j} + dx^2 * dt * F;
    X{j+1} = cqtlyap(M, M, -C, 'debug', true);
end		
	\end{verbatim}
	The parameters \texttt{dt}, \texttt{dx}, and \texttt{timesteps} 
	are set according to the choice of $\Delta t$, $\Delta x$, and the number
	of desired time steps. 
	
		\begin{figure}
		\centering
		\begin{tikzpicture}
		\begin{axis}[view={-40}{40}, width=.5\linewidth,height=.26\textheight, 
		title = {$t = 0.05$}, zmax = 0.05, zmin = 0]
		\addplot3[surf,shader=interp] table {pde_2d_t5.dat};
		\end{axis}
		\end{tikzpicture}~\begin{tikzpicture}
		\begin{axis}[view={-40}{40}, width=.5\linewidth, height=.26\textheight,
		title = {$t = 0.1$}, zmax = 0.05, zmin = 0]
		\addplot3[surf,shader=interp] table {pde_2d_t10.dat};
		\end{axis}
		\end{tikzpicture} \\[5pt]
		\begin{tikzpicture}
		\begin{axis}[view={-40}{40}, width=.5\linewidth, height=.26\textheight,
		title = {$t = 0.15$}, zmax = 0.05, zmin = 0]
		\addplot3[surf,shader=interp] table {pde_2d_t15.dat};
		\end{axis}
		\end{tikzpicture}~\begin{tikzpicture}
		\begin{axis}[view={-40}{40}, width=.5\linewidth, height=.26\textheight,
		title = {$t = 0.2$}, zmax = 0.05, zmin = 0]
		\addplot3[surf,shader=interp] table {pde_2d_t20.dat};
		\end{axis}
		\end{tikzpicture}
		\caption{Plot of the solution restricted to $[0, 2]^2$ 
			for the PDE \eqref{eq:pde} at different
			time steps during the implicit Euler scheme. }
		\label{fig:solpde}
	\end{figure}
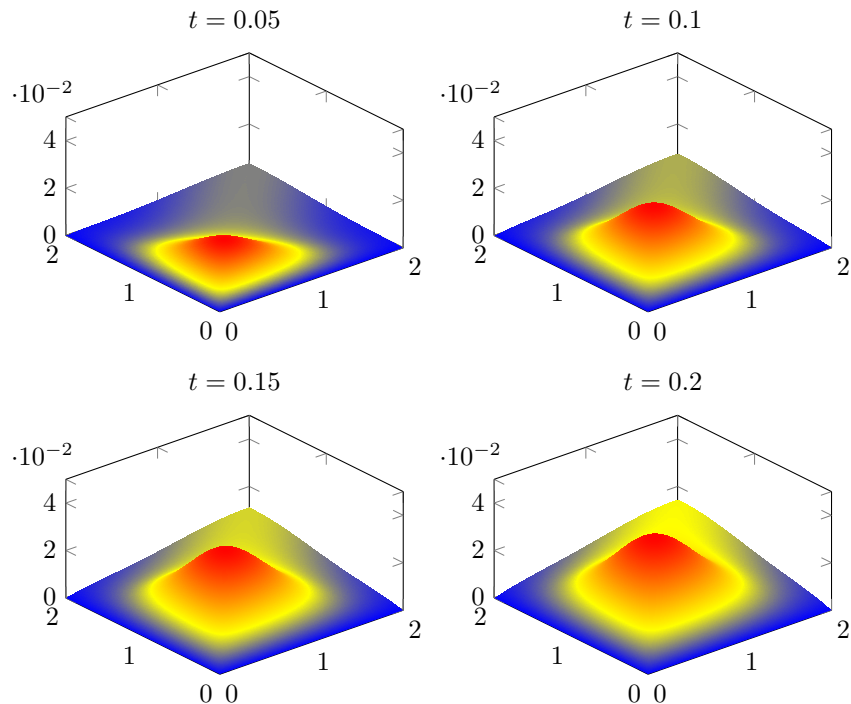
	
	In Figure~\ref{fig:solpde} the part of the solution in $[0, 2] \times [0, 2]$
	is reported at different timesteps; the effects of the term $e^{-(x+y)^2}$ 
	are clearly visible.

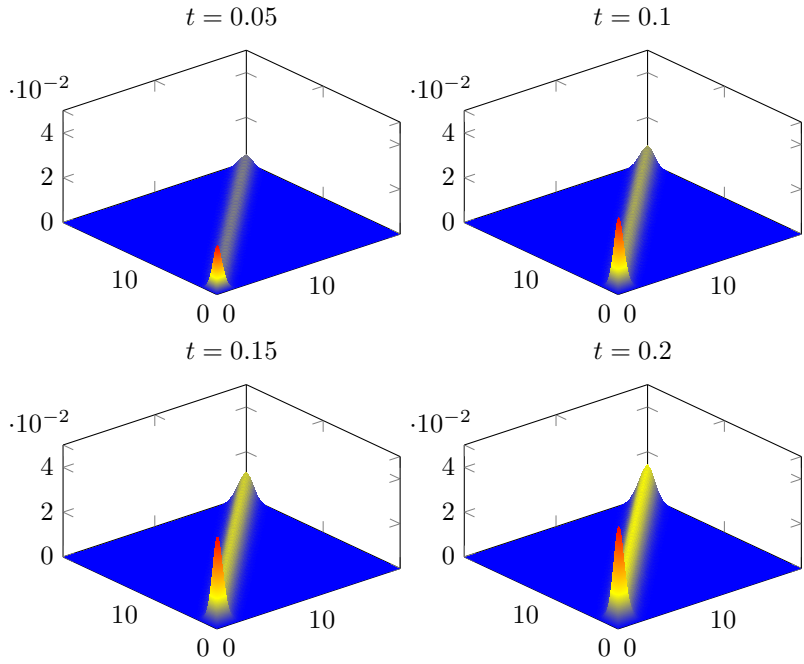
\begin{figure}
	\centering
	\begin{tikzpicture}
	\begin{axis}[view={-40}{40}, width=.5\linewidth, height=.25\textheight,
	title = {$t = 0.05$}, zmax = 0.05, zmin = 0]
	\addplot3[surf,shader=interp] table {pde_2d_t5_large.dat};
	\end{axis}
	\end{tikzpicture}~\begin{tikzpicture}
	\begin{axis}[view={-40}{40}, width=.5\linewidth, height=.25\textheight,
	title = {$t = 0.1$}, zmax = 0.05, zmin = 0]
	\addplot3[surf,shader=interp] table {pde_2d_t10_large.dat};
	\end{axis}
	\end{tikzpicture} \\
	\begin{tikzpicture}
	\begin{axis}[view={-40}{40}, width=.5\linewidth, height=.25\textheight,
	title = {$t = 0.15$}, zmax = 0.05, zmin = 0]
	\addplot3[surf,shader=interp] table {pde_2d_t15_large.dat};
	\end{axis}
	\end{tikzpicture}~\begin{tikzpicture}
	\begin{axis}[view={-40}{40}, width=.5\linewidth, height=.25\textheight,
	title = {$t = 0.2$}, zmax = 0.05, zmin = 0]
	\addplot3[surf,shader=interp] table {pde_2d_t20_large.dat};
	\end{axis}
	\end{tikzpicture}

	\caption{Plot of the solution restricted to $[0, 20]^2$ 
		for the PDE \eqref{eq:pde} at different
		time steps during the implicit Euler scheme.}
	\label{fig:solpdelarge}
\end{figure}

On the other hand, as shown in Figure~\ref{fig:solpdelarge}, 
the solution on
a larger set such as $[0, 20]^2$, is more influenced by 
the Toeplitz part; these plots also suggests that the support of the
solution is unbounded (indeed, the solution $X$ has a nonzero Toeplitz
part, and is therefore not compact). 

\begin{figure}
	\centering
	\begin{tikzpicture}
		\begin{semilogyaxis}[title = {Convergence for the PDE problem},
		  xlabel = Iterations, ylabel = {Relative residual}, width=.75\linewidth,
		  height = .32\textheight]
			\addplot table {pdeconv.dat};
			\addplot table {adipde.dat};
			\addplot[dashed,green] table {pdezolotarev.dat};
			\legend{Galerkin,ADI,ADI error bound~\eqref{eq:adizolotarev}}
		\end{semilogyaxis}	
	\end{tikzpicture}

	\caption{Relative residual in the $2$-norm 
	for the convergence of the Galerkin projection
	method (using poles $0$ and $\infty$), i.e., the extended Krylov solver, 
	and the ADI iteration, using optimal Zolotarev poles. The error
	bound on the relative residual is the one obtained for the ADI 
	method using the error equation, as 
	described in \eqref{eq:adizolotarev}.}
	\label{fig:pdeconv}
\end{figure}
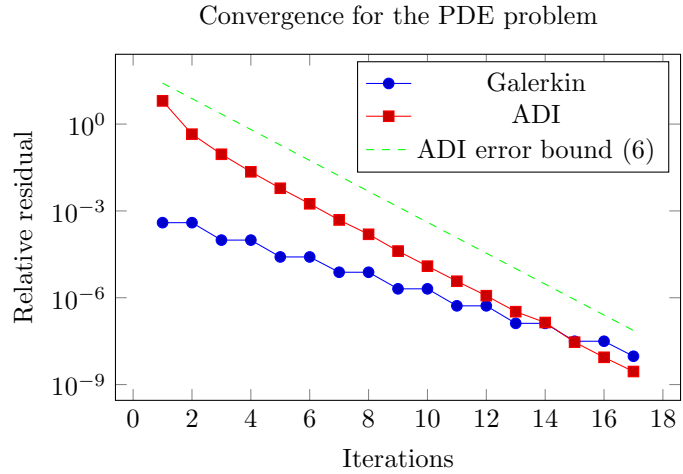

The convergence history for the first step is reported in Figure~\ref{fig:pdeconv}. Note that the Galerkin method approximation
error decreases only every two steps; this is because in the Galerkin 
projection approach we use poles $0$ and $\infty$ alternatively; 
in the method, it is possible to compute the solution and the
residual after adding an infinity pole (see \cite{simoncini2016computational}
for further details), and so adding the pole $0$ gives no immediate benefit, 
but is only useful after the pole $\infty$ is added as well, and a new approximate
solution can be computed. 

We now test the method on a slightly modified example where the exact solution
is known, to check the convergence of the scheme as the discretization steps 
in space and time go to zero. 

We choose $u(x,y,t) = (1 - e^{-t}) (e^{-(x-y)^2} + e^{-x-y})$, which can 
be verified to solve
\begin{equation} \label{eq:pde2}
  \begin{cases}
  \frac{\partial u(x,y,t)}{\partial t} = \Delta u(x,y,t) + f(x,y,t) & 
  \text{on } \mathbb R_+ \times \mathbb R_+ \\
  u(x,y,t) = (1 - e^{-t}) (e^{-(x-y)^2} + e^{-x-y}) & \text{on } \mathbb R_+ \times \{0\} \cup \{ 0\} \times \mathbb R_+ \\
  \end{cases}, 
\end{equation}
where $f(x,y,t)$ is appropriately chosen as follows:
\begin{align*}
	f(x,y,t) &= e^{-x-y} (3e^{-t} - 1) + e^{-(x-y)^2}(e^{-t} + 4(1-e^{-t})(1 - 2(x-y)^2))
\end{align*}
We apply the same approach followed for \eqref{eq:pde}
using different values of $\Delta t = \Delta x$, 
and we check the accuracy compared to the true solution at time 
$t = 1$. We selected $\Delta t = \Delta x$ 
between $\frac{1}{32} = 3.125\cdot 10^{-2}$ and $\frac 12$. 
The results are reported in Table~\ref{fig:pdeconvstep}. 

\begin{table}
	\centering
\begin{tabular}{c|c}
	$\Delta t = \Delta x$ & $\norm{u(1) - \hat u(1)}_\infty$ \\ \hline 
	$0.5$ & $2.4185\cdot 10^{-1}$ \\
	$0.25$ & $1.5782\cdot 10^{-1}$ \\
	$0.125$ & $8.9787\cdot 10^{-2}$ \\
	$6.25 \cdot 10^{-2}$ & $4.9531\cdot 10^{-2}$ \\
	$3.125 \cdot 10^{-2}$ & $2.6063\cdot 10^{-2}$ 
\end{tabular}
	\caption{Error in the infinite norm for the solution of \eqref{eq:pde2} 
	  at the time $t = 1$, using different step sizes in time and 
      space. }
	\label{fig:pdeconvstep}
\end{table}

\subsection{Stein equation for solving quadratic equations}
As another application, we consider the solution of quadratic equations
\[
  A_{-1} + A_0 X + A_1 X^2 = X, 
\]
where $A_j$ are semi-infinite quasi-Toeplitz matrices with non-negative coefficients
and such that $A_{-1} + A_0 + A_1$ is row-stochastic. The natural space where this
problem can be formulated is $\ell^\infty$, and $\qt$ matrices appear when one 
considers quasi-birth-Death problems \cite{bini2005numerical,bini2018quadratic}. Rephrasing the equation as $F(X) = 0$ and using the Newton method yields the
linear equation in $H$:
\[
  (A_0 + A_1 X_k - I) H + A_1 H X_k = -F(X_k)
\]
In this context, the matrix $A_0 + A_1 X_k - I$ is guaranteed to be invertible
\cite{bini2019solving}, and therefore the equation can be recast into a Stein equation
\[
 H + (A_0 + A_1 X_k - I)^{-1} A_1 H X_k = -(A_0 + A_1 X_k - I)^{-1} F(X_k).
\]	
Solving for $H$ gives the next iterate of the Newton method as $X_{k+1} = 
X_k + H$. The
Newton method is often initialized by $X_0 = 0$, since this guarantees convergence. 
Note that with this choice $X_1 = -(A_0 - I)^{-1} A_{-1}$. 
To study the behavior of the different methods in this realistic setting,
we consider
the three different described approaches for the solution of the Stein
equation required to compute the first two Newton iterates, $X_1$ and $X_2$.

We have chosen a test arising from the study of 
a tandem Jackson queue \cite{motyer2006decay}, 
that can be found as model number $3$ in the examples considered in \cite{bini2018quadratic}. This gives rise to slow convergence because it involves
Stein equations with matrices $M,N$ of norms quite close to $1$.

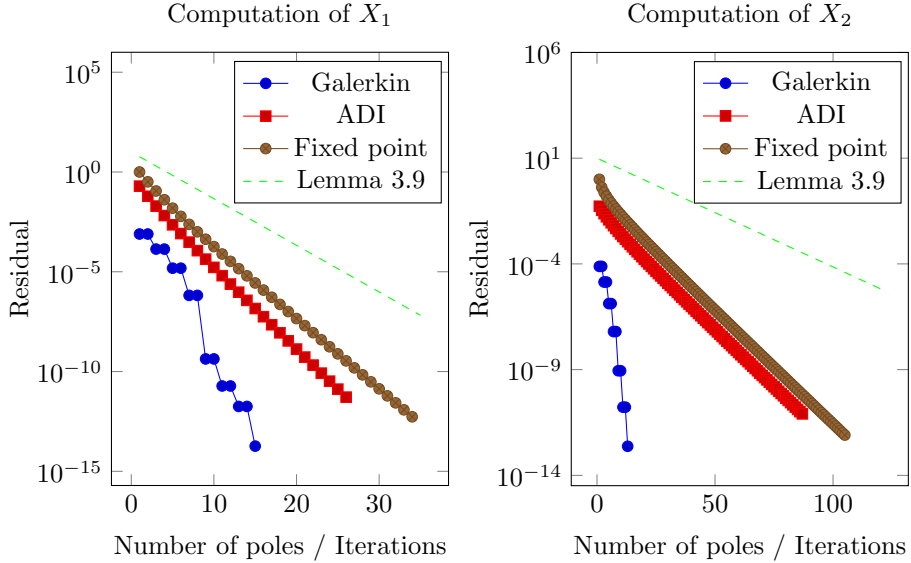
\begin{figure}
	\centering
	\begin{tikzpicture}
	\begin{semilogyaxis}[width=.5\linewidth,
	xlabel = Number of poles / Iterations, ylabel = Residual, ymax = 1000000, 
	height = .38\textheight,
	title = Computation of $X_1$]
	\addplot table {galerkin1.dat};
	\addplot table {adi1.dat};
	\addplot table {fixed1.dat};
	\addplot[dashed,green] table {convX1.dat};
	\legend{Galerkin, ADI, Fixed point,Lemma~\ref{lem:adiconv}}
	\end{semilogyaxis}
	\end{tikzpicture}~\begin{tikzpicture}
	\begin{semilogyaxis}[width=.5\linewidth,
	xlabel = Number of poles / Iterations, ylabel = Residual, ymax = 1000000, 
	height = .38\textheight, title = Computation of $X_2$]
	\addplot table {galerkin.dat};
	\addplot table {adi.dat};
	\addplot table {fixed.dat};
	\addplot[dashed,green] table {convX2.dat};	
	\legend{Galerkin, ADI, Fixed point,Lemma~\ref{lem:adiconv}}
	\end{semilogyaxis}
	\end{tikzpicture}
	
	\caption{Convergence history for the solution of a Stein equation arising
		from the analysis of quasi-Birth-Death processes. On the left it is reported the computation for the first Newton iterate, whereas on the right the
		one for the second one. The convergence rate predicted
		by Lemma~\ref{lem:adiconv} 
		is reported with the dashed green line.}
	\label{fig:stein}
\end{figure}

In particular, one can observe that 
both ADI and ``Fixed point'' (which, as discussed in Section~\ref{sec:stein}, 
are essentially equivalent in this context), require a large number of steps to
converge (especially for $X_2$). 
The Galerkin iteration, on the other hand, has the capability of
adapt to the problem, and manages to converge much faster. 

The disadvantage of the Galerkin method is that the formulation requires
an Hilbert space, and therefore only problem which are well-defined on 
$\ell^2$ are treatable. Since the natural space for the problem arising from
the probabilistic setting is $\ell^\infty$, this might not always be the case.

	In this experiment, the tolerance for all the methods has been set to $10^{-12}$. 
	The implementation can be found in \texttt{cqt-toolbox}, where these equations
	can be solved by using the commands
	\begin{verbatim}
    X = cqtstein(M, N, C, 'method', ...);
	\end{verbatim}
	where \texttt{method} can be either \texttt{'galerkin'}, \texttt{'adi'}, 
	or \texttt{'fixedpoint'}. 
	A refined implementation of the Newton method for 
	the solution of QBD models is beyond the scope of
	this work, and we refer the reader to \cite{bini2019solving} for further
	details.

	\section{Conclusions}
	
	We have discussed the theoretical and numerical tools needed to solve
	linear matrix equation with infinite quasi-Toeplitz matrices. These equations
	appear in a variety of settings, and we have demonstrated that it is
	possible to extend many results known in the finite dimensional
	setting to this more general framework. ADI and fixed point iteration have
	been considered. In addition, we have provided
	an implementation of Galerkin projection scheme for the case when
	the matrices are in $\ell^2$, and we have shown that it can be very effective. 
	
	The methods proposed are made freely available in the MATLAB package
	\texttt{cqt-toolbox}, and provide a promising base for the development
	of Newton methods for quadratic equation that is required in the study
	of Markov chains on infinite bi-dimensional lattices. 
	
	
	\bibliographystyle{elsarticle-harv}
	\bibliography{lyapunov-cqt}
	
\end{document}